**Some results and conjectures about Hankel determinants**

**of sequences which are related to Catalan-like numbers**

Johann Cigler

Fakultät für Mathematik, Universität Wien

johann.cigler@univie.ac.at

**Abstract**

Martin Aigner introduced Catalan-like numbers as elements of the first column of "admissible matrices" and studied Hankel determinants of their forward shifts. In this paper we collect some properties of the Hankel determinants of the other columns which are suggested by computer experiments. By prepending zero rows to admissible matrices we also consider Hankel determinants of backward shifts.

## 1. Introduction and sketch of the situation

Let us first state all results and conjectures and only afterwards provide proofs for the theorems and some hints for the conjectures.

Let $\mathbf{s}=(s_k)_{k\geq 0}$ be a sequence of real numbers. Define numbers $a_{n,k}=a_{n,k}(\mathbf{s})$ by

$$a_{n,k}=a_{n-1,k-1}+s_k a_{n-1,k}+a_{n-1,k+1} \tag{1}$$

with $a_{n,k}=0$ for $k<0$ and $a_{0,k}=[k=0]$.

Following Aigner [1] we call the matrices $\left(a_{n,k}(\mathbf{s})\right)$ *admissible matrices* and the numbers $a_n(\mathbf{s})=a_{n,0}(\mathbf{s})$ *Catalan-like numbers of type* $\mathbf{s}$.

The first numbers are $a_0(\mathbf{s})=1,\ a_1(\mathbf{s})=s_0,\ a_2(\mathbf{s})=s_0^2+1,\ a_3(\mathbf{s})=s_0^3+2s_0+s_1,$ $a_4(\mathbf{s})=s_0^4+3s_0^2+2+2s_0s_1+s_1^2.$

Some well-known special cases are:

$\mathbf{s}=\left(0,0,0,0,\cdots\right)$ gives $a_{2n}(\mathbf{s})=C_n,\ a_{2n+1}(\mathbf{s})=0,$ where $C_n=\dfrac{1}{n+1}\dbinom{2n}{n}$ are the Catalan numbers (OEIS [8], A000108),

$\mathbf{s}=\left(1,0,0,0.\cdots\right)$ gives $a_n(\mathbf{s})=\dbinom{n}{\left\lfloor \frac{n}{2}\right\rfloor},$

$\mathbf{s}=\left(1,2,2,2,\cdots\right)$ gives the Catalan numbers $C_n$,

$\mathbf{s}=\left(2,2,2,2,\cdots\right)$ gives the shifted Catalan numbers $C_{n+1}$,

$\mathbf{s}=\left(1,1,1,1,\cdots\right)$ gives the Motzkin numbers $M_n$ (cf. OEIS [8],A001006).

For a constant type $\mathbf{s}=(c,c,c,\cdots)$ we use the notation $a_{n,k}(c)$ instead of $a_{n,k}(\mathbf{s})$ and for the type $\mathbf{s}=(c+b,c,c,c,\cdots)$ we write $a_{n,k}(b,c)$ and analogous notations for other objects depending on $\mathbf{s}$.

Some aspects for $c=2$ have been considered in [3], [4] and [5]. Here we have $a_{n,k}(2)=[x^n]x^k c(x)^{2k+2}=\frac{k+1}{n+1}\binom{2n+2}{n-k}$ where $c(x)=\frac{1-\sqrt{1-4x}}{2x}$ is the generating function of the Catalan numbers.

We extend $a_{n,k}(\mathbf{s})$ to negative $n$ by setting $a_{n,k}(\mathbf{s})=0$ for $n<0$.

For $m\in\mathbb{Z}$ we define the Hankel matrices

$$\mathbf{A}_{m,k,n}(\mathbf{s})=\left(a_{m+i+j,k}(\mathbf{s})\right)_{i,j=0}^{n-1} \tag{2}$$

and the Hankel determinants

$$D_{m,k,n}(\mathbf{s})=\det\mathbf{A}_{m,k,n}(\mathbf{s}). \tag{3}$$

For $k=0$ we also write $\mathbf{A}_{m,0,n}(\mathbf{s})=\mathbf{A}_{m,n}(\mathbf{s})$ and $D_{m,0,n}(\mathbf{s})=D_{m,n}(\mathbf{s})$.

Let

$$\mathbf{D}_m(\mathbf{s})=\left(D_{m,0,n}(\mathbf{s})\right)_{n\geq 0} \tag{4}$$

and for $\lambda\in\mathbb{R}$ and a positive integer $r$

$$\lambda T_r\left(\mathbf{D}_m(\mathbf{s})\right):=T_r\left(\lambda\mathbf{D}_m(\mathbf{s})\right)=\left(1,0,\cdots,0,\lambda D_{m,0}(\mathbf{s}),\lambda D_{m,1}(\mathbf{s}),\lambda D_{m,2}(\mathbf{s}),\cdots\right), \tag{5}$$

i.e. $\lambda T_r\left(\mathbf{D}_m(\mathbf{s})\right)=\left(u_n\right)_{n\geq 0}$ with $u_0=1,\ u_n=0$ for $0<n<r$ and $u_{r+i}=\lambda D_{m,i}(\mathbf{s})$ for $i\in\mathbb{N}$.

For Catalan-like numbers the Hankel determinants $\mathbf{D}_{-m}(\mathbf{s})$ are given by

**Theorem 1**

*For* $m\in\mathbb{N}$ *we have*

$$\mathbf{D}_{-m}(\mathbf{s})=\left(D_{-m,0,n}(\mathbf{s})\right)_{n\geq 0}=(-1)^{\binom{m+1}{2}}T_{m+1}\left(\mathbf{D}_m(E\mathbf{s})\right). \tag{6}$$

*with* $E\mathbf{s}=\left(s_1,s_2,s_3,\cdots\right)$.

**Remark**

For $m=0$ this is true because $\mathbf{D}_0(\mathbf{s})=(1,1,1,\cdots)$.

For example, for $\mathbf{s}=(1,0,0,0,\cdots)$ we have $\left(D_{2,0,n}(1,0)\right)_{n\geq 0}=(1,2,3,4,5,6,7,8,9,10,11,12,\cdots)$,

for $E\mathbf{s}=(0,0,0,\cdots)$ we get $\left(D_{2,0,n}(0)\right)_{n\geq 0}=(1,1,2,2,3,3,4,4,5,5,6,6,\cdots)$, and

$$\left(D_{-2,0,n}(1,0)\right)_{n\geq 0}=\left(1,0,0,-1,-1,-2,-2,-3,-3,-4,-4,-5,-5,-6,-6,-7,-7,-8,-8,\cdots\right).$$

For constant types there exists an analog of Theorem 1. Let

$$\mathbf{D}_{m,k}(c)=\left(D_{m,k,n}(c)\right)_{n\geq 0}. \tag{7}$$

**Conjecture 2**

*For* $m\in\mathbb{N}$ *we have*

$$\mathbf{D}_{-m,k}(c)=(-1)^{\binom{m+k+1}{2}}T_{m+k+1}\left(\mathbf{D}_{m,k}(c)\right). \tag{8}$$

For example,

$$\mathbf{D}_{2,3}(0)=\left(1,0,-1,0,4,0,-4,0,9,0,-9,0,16,0,-16,0,25,0,-25,\cdots\right),$$

$$\mathbf{D}_{-2,3}(0)=-T_6\left(\mathbf{D}_{2,3}(0)\right)=\left(1,0,0,0,0,0,-1,0,1,0,-4,0,4,0,-9,0,9,0,-16,\cdots\right).$$

For $m=0$ formula (8) seems to be the real reason for

**Theorem 3 ([6])**

*For* $k\geq 1$

$$\begin{aligned} &D_{0,k,(k+1)n}(c)=(-1)^{n\binom{k+1}{2}}, \\ &D_{0,k,n}(c)=0 \quad \text{else.} \end{aligned} \tag{9}$$

For assume Conjecture 2 and let $\mathbf{D}_{0,k}(c)=\left(d_n\right)_{n\geq 0}$. Then (8) gives $d_0=1$ and $d_n=0$ for $0<n<k+1$ and $d_{n+k+1}=(-1)^{\binom{k+1}{2}}d_n$.

It turns out that all Hankel determinants $D_{m,k,n}(c)$ show a modular pattern modulo $k+1$.

For $m=1$ this can be seen from

**Theorem 4 ([6])**

*For* $k\geq 1$

$$\begin{aligned} &D_{1,k,(k+1)n}(c)=(-1)^{n\binom{k+1}{2}}F_{n+1}\left(L_{k+1}(c)\right), \\ &D_{1,k,(k+1)n+k}(c)=(-1)^{n\binom{k+1}{2}+\binom{k}{2}}F_{n+1}\left(L_{k+1}(c)\right), \\ &D_{1,k,n}(c)=0 \quad \text{if } n\not\equiv 0,k \bmod (k+1). \end{aligned} \tag{10}$$

The original proof, which we show later, expresses the determinants in terms of $x,y$ satisfying $x+y=c$ and $xy=1$. But in the present context we find it more convenient to formulate them with variants of the Fibonacci polynomials $F_n(x)$ and Lucas polynomials $L_n(x)$ defined by

$F_n(x) = xF_{n-1}(x) - F_{n-2}(x)$ with initial values $F_0(x) = 0$ and $F_1(x) = 1$ and
$L_n(x) = xL_{n-1}(x) - L_{n-2}(x)$ with initial values $L_0(x) = 2$ and $L_1(x) = x.$

For $m = 2$ we get

**Conjecture 5**

*For* $k \geq 1$ *the Hankel determinants* $D_{2,k,n}(c)$ *are*

$$\begin{aligned} D_{2,k,(k+1)n}(c) &= (-1)^{n\binom{k+1}{2}} F_{n+1}\left(L_{k+1}(c)\right)^2, \\ D_{2,k,(k+1)n+k-1}(c) &= (-1)^{n\binom{k+1}{2}+\binom{k-1}{2}} F_{n+1}\left(L_{k+1}(c)\right)^2, \end{aligned} \tag{11}$$

$$D_{2,k,(k+1)n+k}(c) = (-1)^{n\binom{k+1}{2}+\binom{k}{2}} (k+1)F_{k+1}(c) \sum_{j=0}^{n} F_{j+1}\left(L_{k+1}(c)\right)^2, \tag{12}$$

*and* $D_{2,k,n}(c) = 0$ *else.*

**Remark**

Formula (12) is a natural analog of the well-known formula

$$D_{2,0,n}(c) = \sum_{j=0}^{n} F_{j+1}(c)^2. \tag{13}$$

The first lines of (9), (10), (11) are special cases of

**Conjecture 6**

*For* $0 \leq m \leq k+1$

$$D_{m,k,(k+1)n}(c) = (-1)^{n\binom{k+1}{2}} F_{n+1}\left(L_{k+1}(c)\right)^m. \tag{14}$$

**Remark**

For $c \in \{0,1,2\}$ the sequence $\left(L_n(c)\right)_{n\geq 0}$ is periodic.

$\left(L_n(0)\right)_{n\geq 0} = \left(\underline{2}, 0, -2, \underline{0}, \cdots\right)$, $\left(L_n(1)\right)_{n\geq 0} = \left(\underline{2}, 1, -1, -2, -1, \underline{1}, \cdots\right)$, $\left(L_n(2)\right)_{n\geq 0} = \left(\underline{2}, \cdots\right)$.

This gives some simplifications in the above formulas. For example, for $c = 2$ $F_{n+1}\left(L_{k+1}(2)\right) = n+1$ and (14) becomes

$$D_{m,k,(k+1)n}(2) = (-1)^{n\binom{k+1}{2}} (n+1)^m \tag{15}$$

for $0 \leq m \leq k+1$.

Finally, we consider types of the form $\mathbf{s} = (c+b, c, c, c, \cdots)$.

It is easy to see that

$$\begin{aligned} D_{0,0,n}(b,c) &= 1, \\ D_{1,0,n}(b,c) &= F_{n+1}(c) + bF_n(c), \\ D_{2,0,n}(b,c) &= \sum_{j=0}^{n} \left(F_{j+1}(c) + bF_j(c)\right)^2. \end{aligned} \tag{16}$$

For $k > 0$ and general $b$ we did not find nice formulas. Whereas $D_{0,k,n}(c)$ takes only values $0$ and $\pm 1$ independent of $c$ it seems that $D_{0,k,n}(b,c)$ does not depend on $c$ either, but is in general a nonvanishing polynomial in $b$.

For example

$$D_{0,1,n}(b,c) = \left(1, 0, -1, b, 1-b^2, -2b+b^3, -1+3b^2-b^4, 3b-4b^3+b^5, \cdots\right).$$

$$D_{0,2,n}(b,c) = \left(1, 0, 0, -1, 0, b^2, 1-b^2, b^2-b^4, -3b^2+2b^4, -1+3b^2-3b^4+b^6, \cdots\right).$$

For $k = 1$ we get $D_{0,1,n}(b,c) = -F_n(-b)$. It seems that the sequences $h_{n,k} = D_{0,k,n}(b,c)$ satisfy recurrence relations of order $2^k$. The first recurrences are

$$\begin{aligned} & h_{n,1} + bh_{n+1,1} + h_{n+2,1} = 0, \\ & h_{n,2} + h_{n+1,2} + b^2 h_{n+2,2} + h_{n+3,2} + h_{n+4,2} = 0, \\ & h_{n,3} - bh_{n+1,3} - \left(b^3 - b\right) h_{n+3,3} + \left(b^4 + b^2 - 2\right) h_{n+4,3} - \left(b^3 - b\right) h_{n+5,3} - bh_{n+7,3} + h_{n+8,3} = 0. \end{aligned} \tag{17}$$

They also seem to satisfy

$$h_{n,k} h_{n-2k-2,k} - b^2 h_{n-1,k} h_{n-1-2k,k} = (1-b^2) h_{n-1-k}^2. \tag{18}$$

For $k = 1$ this reduces to $F_n(x) F_{n-4}(x) - x^2 F_{n-1}(x) F_{n-3}(x) = \left(1 - x^2\right) F_{n-2}(x)^2$.

This follows from Cassini's formula $F_{n-2}(x)^2 - F_{n-1}(x) F_{n-3}(x) = 1$ and the easily verified formula $F_n(x)^2 - F_n(x) F_{n-4}(x) = x^2$.

For $b = 0$ and $b = \pm 1$ the sequences $\left(D_{0,k,n}(b,c)\right)_{n \ge 0}$ are periodic. For $b = 0$ with period $2(k+1)$ and for $b = \pm 1$ with period $2k+1$. This explains the modular pattern of the determinants $\left(D_{m,k,n}(b,c)\right)_{n \ge 0}$ for these cases. Since $a_{n,k}(-\mathbf{s}) = (-1)^{n-k} a_{n,k}(\mathbf{s})$ it suffices in the following to consider $b = 1$.

**Conjecture 7**

$$\begin{aligned} D_{0,k,(2k+1)n}(1,c) &= 1 \\ D_{0,k,(2k+1)n+k+1}(1,c) &= (-1)^{\binom{k+1}{2}}. \end{aligned} \tag{19}$$

and $D_{0,k,n}(1,c)=0$ else.

For the determinants $D_{1,k,n}(1,c)$ we get

**Conjecture 8**

*For* $k>0$

$$D_{1,k,(2k+1)n}(1,c)=F_n\left(L_{2k+1}(c)\right)+F_{n+1}\left(L_{2k+1}(c)\right) \tag{20}$$

$$D_{1,k,(2k+1)n+k}(1,c)=(-1)^{\binom{k}{2}}d_{1,k,(2k+1)n}(1,c) \tag{21}$$

$$D_{1,k,(2k+1)n+k+1}(1,c)=(-1)^{\binom{k+1}{2}}\left(L_k(c)+L_{k+1}(c)\right)F_{n+1}\left(L_{2k+1}(c)\right), \tag{22}$$

$$D_{1,k,(2k+1)n+2k}(1,c)=(-1)^k\left(L_k(c)+L_{k+1}(c)\right)F_{n+1}\left(L_{2k+1}(c)\right). \tag{23}$$

*and* $D_{1,k,n}(1,c)=0$ *else.*

**Conjecture 9**

*For* $k\geq 3$ *the determinants* $D_{2,k,(2k+1)n+r}(1,c),\ 0\leq r\leq 2k,$ *vanish except for* $r=0,k-1,k,k+1,2k-1,2k.$

*For* $k=1$ *and* $c\geq 0$ *we get*

$$D_{2,1,3n}(1,c)=\left(F_n(L_3(c))+F_{n+1}(L_3(c))\right)^2 \tag{24}$$

$$D_{2,1,3n+2}(1,c)+D_{2,1,3n+1}(1,c)=-F_{n+1}\left(L_3(c)\right)^2(c+2)^2\left(F_2(c)-F_1(c)\right)^2 \tag{25}$$

*For* $k\geq 2$ *and* $c\geq 1$ *we get*

$$D_{2,k,(2k+1)n}(1,c)=\left(F_n(L_{2k+1}(c))+F_{n+1}(L_{2k+1}(c))\right)^2, \tag{26}$$

$$D_{2,k,(2k+1)n+k-1}(1,c)=(-1)^{\binom{k-1}{2}}\left(F_n(L_{2k+1}(c))+F_{n+1}(L_{2k+1}(c))\right)^2, \tag{27}$$

$$D_{2,k,(2k+1)n+k+1}(1,c)=(-1)^{\binom{k+1}{2}}F_{n+1}\left(L_{2k+1}(c)\right)^2(c+2)^2\left(F_{k+1}(c)-F_k(c)\right)^2 \tag{28}$$

$$D_{2,k,(2k+1)n+2k-1}(1,c)=-F_{n+1}\left(L_{2k+1}(c)\right)^2(c+2)^2\left(F_{k+1}(c)-F_k(c)\right)^2 \tag{29}$$

$$\begin{aligned}&D_{2,k,(2k+1)n+2k}(1,c)+(-1)^{\binom{k-1}{2}}D_{2,k,(2k+1)n+k}(1,c)\\&=(-1)^k\left(F_{n+1}\left(L_{2k+1}(c)\right)\right)^2(c+2)^2\left(F_{k+1}(c)-F_k(c)\right)^2\left(F'_{k+1}(c)-F'_k(c)\right).\end{aligned} \tag{30}$$

It seems that there exists no analog of Conjecture 2 for $D_{-m,k,n}(b,c)$ except for $D_{m,k,n}(-1,2)$ where

$$D_{1-m,k,n+m+k}(-1,2)=(-1)^{\binom{m+k}{2}}D_{m,k,n}(-1,2). \tag{31}$$

Only for $m=1$ we found the following identities which could be interpreted as such an analog:

$$\begin{aligned} D_{-1,k,(2k+1)(n+1)}(1,c) &= -D_{1,k,(2k+1)n}(1,c) \\ D_{-1,k,(2k+1)n+k+2}(1,c) &= -D_{1,k,(2k+1)n+k}(1,c) \end{aligned} \tag{32}$$

$$\begin{aligned} D_{-1,k,(2k+1)n+1}(1,c) &= -\frac{D_{1,k,(2k+1)n+2k}(1,c)}{c+2}, \\ D_{-1,k,(2k+1)n+k+1}(1,c) &= -\frac{D_{1,k,(2k+1)n+k+1}(1,c)}{c+2}. \end{aligned} \tag{33}$$

All other determinants vanish.

For example, for $k=1$ and $c=3$ we get

$$\left(D_{1,1,n}(1,3)\right)_{n\geq 0}=\left(1,1,-10,19,19,-180,341,341,-3230,6119,6119,-57960,109801,\cdots\right)$$

$$\left(D_{-1,1,n}(1,3)\right)_{n\geq 0}=\left(1,0,0,-1,2,2,-19,36,36,-341,646,646,-6119,11592,11592,-109801,208010,\cdots\right).$$

$$\left(D_{1,1,3n}(1,3)\right)_{n\geq 0}=\left(D_{1,1,3n+1}(1,3)\right)_{n\geq 0}=\left(1,19,341,6119,109801,1970299,35355581,\cdots\right),$$
$$\left(D_{-1,1,3(n+1)}(1,3)\right)_{n\geq 0}=\left(-1,-19,-341,-6119,-109801,-1970299,-35355581,\cdots\right),$$

$$\left(D_{1,1,3n+2}(1,3)\right)_{n\geq 0}=\left(-10,-180,-3230,-57960,-1040050,-18662940,-334892870,\cdots\right)$$
$$=5\left(-2,-36,-646,-11592,-208010,-3732588,-66978574,\cdots\right),$$

$$\left(D_{-1,1,3(n+1)+1}(1,3)\right)_{n\geq 0}=\left(D_{-1,1,3(n+1)+2}(1,3)\right)_{n\geq 0}=\left(2,36,646,11592,208010,3732588,66978574,\cdots\right).$$

## 2. Some well-known facts

The generating function of $\left(a_n(\mathbf{s})\right)$ is (cf. [2], Theorem 5.8.2)

$$A(x,\mathbf{s})=\sum_{n\geq 0}a_n(\mathbf{s})x^n=\frac{1}{1-s_0x-x^2A(x,E\mathbf{s})} \tag{34}$$

where $E\mathbf{s}=\left(s_1,s_2,\cdots\right)$.

As shown in [1] we have

$$D_{0,n}(\mathbf{s})=1,$$
$$D_{1,n}(\mathbf{s})=s_{n-1}D_{1,n-1}-D_{1,n-2}. \tag{35}$$

Let us recall the proofs.

We have

$$\sum_{k=0}^{\min(m,n)} a_{n,k}a_{m,k}=a_{n+m,0}. \tag{36}$$

For (36) is trivially true for $n=0$ and arbitrary *m*. If it is true for $n-1$ and arbitrary $m$ then we get

$$\sum_k a_{n,k}a_{m,k}=\sum_k a_{m,k}\left(a_{n-1,k-1}+s_k a_{n-1,k}+a_{n-1,k+1}\right)$$
$$=\sum_k a_{n-1,k}\left(a_{m,k+1}+s_k a_{m,k}+a_{m,k-1}\right)$$
$$=\sum_k a_{n-1,k}a_{m+1,k}=a_{m+n,0}.$$

Identity (36) is equivalent with

$$A_nA_n^T=\left(a_{i+j,0}\right)_{i,j=0}^{n-1} \tag{37}$$

where $A_n=\left(a_{i,j}\right)_{i,j=0}^{n-1}$.

Since $A_n$ is triangular with $a_{n,n}=1$ for $n\in\mathbb{N}$ we get

$$\det\left(a_{i+j,0}\right)_{i,j=0}^{n-1}=1. \tag{38}$$

Let $B_n=\left(a_{i+1,j}\right)_{i,j=0}^{n-1}$. Then by (36) $B_nA_n^T=\left(a_{i+j+1,0}\right)_{i,j=0}^{n-1}$. Therefore $D_{1,n}(\mathbf{s})=\det B_n$.

Let

$$J_n=\begin{pmatrix} s_0 & 1 & 0 & \cdots & 0\\ 1 & s_1 & 1 & \cdots & 0\\ 0 & 1 & s_2 & \cdots & 0\\ & & \cdots & & \\ 0 & 0 & 0 & \cdots & s_{n-1}\end{pmatrix} \tag{39}$$

be a tridiagonal matrix with $s_i$ in the main diagonal and 1 in the adjacent diagonals.

Then

$$B_n=A_nJ_n \tag{40}$$

because

$$A_nJ_n = \left(a_{i,j-1} + s_j a_{i,j} + a_{i,j+1}\right)_{i,j=0}^{n-1} = \left(a_{i+1,j}\right)_{i,j=0}^{n-1} = B_n.$$

This implies

$$D_{1,n}(\mathbf{s}) = \det B_n = \det J_n. \tag{41}$$

Expanding $\det J_n$ with respect to the last column gives

$$D_{1,n}(\mathbf{s}) = s_{n-1}D_{1,n-1}(\mathbf{s}) - D_{1,n-2}(\mathbf{s}) \tag{42}$$

with initial values $D_{1,0}(\mathbf{s}) = 1$ and $D_{1,1}(\mathbf{s}) = s_0$.

### 3. Proof of Theorem 1

We want to show that

$$\begin{aligned} &D_{-m,n}(\mathbf{s}) = 0 \ \text{ for } 0 < n \le m, \\ &D_{-m,n+m+1}(\mathbf{s}) = (-1)^{\binom{m+1}{2}} D_{m,n}(E\mathbf{s}) \ \text{ for } n \in \mathbb{N}. \end{aligned} \tag{43}$$

$\det \mathbf{A}_{-m,n}(\mathbf{s}) = 0$ for $1 \le n \le m$, because all elements $a_{-m+j}(\mathbf{s}),\ 0 \le j \le n-1,$ of the top row vanish and

$$\det \mathbf{A}_{-m,m+1}(\mathbf{s}) = (-1)^{\binom{m+1}{2}}, \tag{44}$$

because $\mathbf{A}_{-m,m+1}(\mathbf{s})$ is a right-triangular matrix with all diagonal entries $=1$.

For example, $\mathbf{A}_{-2,3}(\mathbf{s}) = \begin{pmatrix} 0 & 0 & -1 \\ 0 & -1 & a_1(\mathbf{s}) \\ -1 & a_1(\mathbf{s}) & a_2(\mathbf{s}) \end{pmatrix} = \begin{pmatrix} 0 & 0 & -1 \\ 0 & -1 & s_0 \\ -1 & s_0 & s_0^2+1 \end{pmatrix}.$

We set $\mathbf{V}_{k,n} = \mathbf{A}_{k-n,n}(\mathbf{s})$ and $v_{k,n} = \det \mathbf{V}_{k,n}$.

Note that the top row of $\mathbf{V}_{k,n} = \mathbf{A}_{k-n,n}(\mathbf{s})$ has precisely $k$ non-vanishing entries $a_{k-n+j}(\mathbf{s}),$ $n-k \le j \le n-1$. For example, $\mathbf{A}_{-2,3}(\mathbf{s}) = \mathbf{V}_{1,3}$.

We know that

$$v_{1,m+1} = (-1)^{\binom{m+1}{2}}. \tag{45}$$

Next we show that

$$v_{2,m+2} = D_{-m,m+2}(\mathbf{s}) = (-1)^{\binom{m+1}{2}} a_{m+1}(E\mathbf{s}). \tag{46}$$

Observe that $\mathbf{V}_{2,n}$ can be obtained from $\mathbf{V}_{1,n+1}$ by deleting the first row and column. By Cramer's rule

$\mathbf{V}_{1,n+1}^{-1} = \frac{1}{\det \mathbf{V}_{1,n+1}} \left(\alpha_{j,i}\right)_{i,j=0}^{n}$ with $\alpha_{j,i} = (-1)^{i+j} \det A_{j,i}$, where $A_{i,j}$ is the matrix obtained by crossing out row $i$ and column $j$ in $\mathbf{V}_{1,n+1}$. Thus $A_{0,0} = \mathbf{V}_{2,n}$.

Therefore $v_2(n)$ is the entry in position $(0,0)$ of $\mathbf{V}_{1,n+1}^{-1}$.

To obtain the inverse we make use of

**Lemma 10**

*Let* $a(n) = b(n) = c(n) = 0$ *for* $n < 0$. *Then*

$\left(a(i+j-n)\right)_{i,j=0}^{n} \left(b(n-j-k)\right)_{j,k=0}^{n} = \left(c(i-k)\right)_{i,k=0}^{n}$ *if and only if*
$\sum_{n\geq 0} a(n)x^n \sum_{n\geq 0} b(n)x^n = \sum_{n\geq 0} c(n)x^n.$

For example

$$\begin{pmatrix} 0 & 0 & 0 & a(0) \\ 0 & 0 & a(0) & a(1) \\ 0 & a(0) & a(1) & a(2) \\ a(0) & a(1) & a(2) & a(3) \end{pmatrix} \cdot \begin{pmatrix} b(3) & b(2) & b(1) & b(0) \\ b(2) & b(1) & b(0) & 0 \\ b(1) & b(0) & 0 & 0 \\ b(0) & 0 & 0 & 0 \end{pmatrix} = \begin{pmatrix} c(0) & 0 & 0 & 0 \\ c(1) & c(0) & 0 & 0 \\ c(2) & c(1) & c(0) & 0 \\ c(3) & c(2) & c(1) & c(0) \end{pmatrix}.$$

The proof is obvious because
$\sum a(i+j-n)b(n-j-k) = \sum_{j\geq n-i} a(j+i-n)b(n-j-k) = \sum_{j\geq 0} a(j)b(i-k-j) = c(i-k).$

Therefore, the inverse of $\mathbf{V}_{1,n+1} = \left(a_{i+j-n}(s)\right)_{i,j=0}^{n}$ is $\mathbf{V}_{1,n+1}^{-1} = \left(b_{n-i-j}\right)_{i.j=0}^{n}$ where

$\sum_{n\geq 0} b_n x^n = \frac{1}{A(x,\mathbf{s})} = 1 - s_0 x - x^2 A(x, E\mathbf{s}) = 1 - s_0 x - \sum_{n\geq 2} a_{n-2}(E\mathbf{s})x^n.$

Thus for $n \geq 2$ we get $b_n = -a_{n-2}(E\mathbf{s})$ and therefore $v_{2,n} = (-1)^{\binom{n-1}{2}} a_{n-2}(E\mathbf{s}).$

For example,

$$v_{2,3} = \det \begin{pmatrix} 0 & 1 & s_0 \\ 1 & s_0 & 1+s_0^2 \\ s_0 & 1+s_0^2 & 2s_0+s_0^3+s_1 \end{pmatrix} = \det \begin{pmatrix} 0 & 1 & 0 \\ 1 & s_0 & 1 \\ s_0 & 1+s_0^2 & s_0+s_1 \end{pmatrix} = -s_1.$$

The condensation formula (cf.[7], Prop. 10)

$$\det A \det A_{1,n}^{1,n} = \det A_1^1 \det A_n^n - \det A_1^n \det A_n^1 \qquad (47)$$

where $A_{i_1,i_2,\cdots,i_k}^{j_1,j_2,\cdots,j_k}$ denotes the submatrix of $A$ in which rows $i_1,i_2,\cdots,i_k$ and columns $j_1,j_2,\cdots,j_k$ are omitted, gives

$$v_k(n+k)v_k(n+k-2)=v_{k-1}(n+k-1)v_{k+1}(n+k-1)-v_k(n+k-1)^2 \tag{48}$$

We want to prove that $v_k(n)=(-1)^{\binom{n-k+1}{2}}D_{n-k,k-1}\left(E\mathbf{s}\right)$.

For the proof let us assume that $s_0,s_1,s_2,\cdots$ are indeterminates and set

$v_k(n)=(-1)^{\binom{n-k+1}{2}}d_{n-k,k-1}$. Then (48) is equivalent with

$-d_{n,k-1}d_{n-2,k-1}+d_{n,k-2}d_{n-2,k}+d_{n-1,k-1}^2=0$. Each $d_{n,k}$ is a polynomial in $s_1,s_2,s_3,\cdots$.

On the other hand by condensation we also know that

$-D_{n,k-1}(E\mathbf{s})D_{n-2,k-1}(E\mathbf{s})+D_{n,k-2}(E\mathbf{s})D_{n-2,k}(E\mathbf{s})+D_{n-1,k-1}^2(E\mathbf{s})=0.$

Since for $\mathbf{s}=\left(2,2,2,\cdots\right)$ all $D_{n,k}(\mathbf{s})\neq 0$ no $d_{n,k}$ vanishes. Therefore we get

$$v_k(n)=\frac{v_{k-1}(n+1)v_{k-1}(n-1)+v_{k-1}(n)^2}{v_{k-2}(n)} \tag{49}$$

for $n\geq k$ with $v_k(k)=D_{0,k}(\mathbf{s})=1$.

Since $v_1(n)=(-1)^{\binom{n}{2}}$ and $v_2(n)=(-1)^{\binom{n-1}{2}}a_{n-2}(E\mathbf{s})$ all $v_k(n)$ for $n\geq k$ are uniquely determined and we get $v_k(n)=D_{k-n,n}(\mathbf{s})=(-1)^{\binom{n-k+1}{2}}D_{n-k,k-1}\left(E\mathbf{s}\right)$.

## 4. The special case $\mathbf{s}=\left(c,c,c,\cdots\right)$.

For $\mathbf{s}=\left(c,c,c,\cdots\right)$ we replace $\mathbf{s}$ with $c$ in all formulas.

By (34) we get $A(x,c)=\dfrac{1}{1-cx-x^2A(x,c)}$ which gives

$$A(x,c)=\frac{1-cx-\sqrt{1-2cx+(c^2-4)x^2}}{2x^2} \tag{50}$$

and

$$A(x,c)=1+cxA(x,c)+x^2A(x,c)^2. \tag{51}$$

Then $\left(a_n(c)\right)_{n\geq 0}=\left(1,c,c^2+1,c^3+3c,c^4+6c^2+2,\cdots\right)$.

The generating function of $\left(a_{n,k}(c)\right)_{n\geq 0}$ is

$$\sum_{n\geq 0} a_{n,k}(c)x^n = x^k A(x,c)^{k+1}. \tag{52}$$

For by (51) we get

$$x^k A(x,c)^{k+1} = x\left(x^{k-1}A(x,c)^k + c\left(x^k A(x,c)^{k+1}\right) + \left(x^{k+1}A(x,c)^{k+2}\right)\right).$$

Comparing coefficients with (1) gives (52).

For $k=0$ we have $D_{0,n}(c)=1,$ by (42) $D_{1,n}(c) = cD_{1,n-1}(c) - D_{1,n-2}(c)$ with initial values $D_{1,0}(c)=1,\ D_{1,1}(c)=c.$

Therefore

$$\begin{aligned} D_{0,n}(c) &= 1, \\ D_{1,n}(c) &= F_{n+1}(c). \end{aligned} \tag{53}$$

From the condensation formula (47) we get

$$D_{2,n}(c)D_{0.n}(c) - D_{0,n+1}(c)D_{2,n-1}(c) - D_{1,n}(c)^2 = 0.$$

By (53) this gives $D_{2,n}(c) - D_{2,n-1}(c) - F_{n+1}(c)^2 = 0$ and therefore

$$D_{2,n}(c) = \sum_{j=0}^{n+1} F_j(c)^2. \tag{54}$$

Finally let us recall the proof of Theorem 4.

In [6], Theorem 1 and Theorem 2 the following result has been proved:

Define numbers $c_{n,k} = c_{n,k}(x,y)$ by

$$c_{n,k} = c_{n-1,k-1} + (x+y)c_{n-1,k} + xyc_{n-1,k+1} \tag{55}$$

with $c_{n,k}=0$ for $k<0$ and $c_{0,k}=[k=0].$

Then

$$\begin{aligned} \det\left(c_{i+j,k}(x,y)\right)_{i,j=0}^{n-1} &= (-1)^{m\binom{k+1}{2}}(xy)^{(k+1)^2\binom{m}{2}} \text{ for } n=(k+1)m, \\ \det\left(c_{i+j,k}(x,y)\right)_{i,j=0}^{n-1} &= 0 \quad \text{for } n \not\equiv 0 \bmod (k+1). \end{aligned} \tag{56}$$

and

$$\det\left(c_{i+j+1,k}(x,y)\right)_{i,j=0}^{(k+1)m-1} = (-1)^{m\binom{k+1}{2}}(xy)^{(k+1)^2\binom{m}{2}}\frac{y^{(k+1)(m+1)}-x^{(k+1)(m+1)}}{y^{k+1}-x^{k+1}},$$
$$\det\left(c_{i+j+1,k}(x,y)\right)_{i,j=0}^{(k+1)m+k-1} = (-1)^{m\binom{k+1}{2}+\binom{k}{2}}(xy)^{(k+1)^2\binom{m}{2}+mk(k+1)}\frac{y^{(k+1)(m+1)}-x^{(k+1)(m+1)}}{y^{k+1}-x^{k+1}} \tag{57}$$
$$\det\left(c_{i+j,k}(x,y)\right)_{i,j=0}^{n-1} = 0 \quad \text{for } n \not\equiv 0, k \bmod (k+1).$$

If we choose $x = \dfrac{c-\sqrt{c^2-4}}{2}$, $y = \dfrac{c+\sqrt{c^2-4}}{2}$, then $x+y=c$ and $xy=1$.

Therefore we get from (56)

$$D_{0,k,(k+1)n}(c) = (-1)^{n\binom{k+1}{2}},$$
$$D_{0,k,n}(c) = 0 \quad \text{else.} \tag{58}$$

Observe that

$$F_n(c) = \frac{y^n - x^n}{y-x} \tag{59}$$

and

$$L_n(c) = x^n + y^n. \tag{60}$$

To see this note that $(z-x)(z-y) = z^2 - cz + 1$ and therefore $x^n = cx^{n-1} + x^{n-2}$ and $y^n = cy^{n-1} + y^{n-2}$.

For $xy=1$ we get $y^{n+2} - x^{n+2} = (x+y)\left(y^{n+1} - x^{n+1}\right) - \left(y^n - x^n\right)$.

Therefore changing $x \to x^k$ and $y \to y^k$ we get

$$\frac{y^{kn}-x^{kn}}{y^k-x^k} = F_n(x^k+y^k) = F_n\left(L_k(c)\right).$$

Therefore (57) implies (10).

## References

[1] Martin Aigner, Catalan-like numbers and determinants, J. Comb. Th. A 87 (1999), 33-51

[2] George E. Andrews, Richard Askey and Ranjan Roy, Special Functions, Encyclopedia of Mathematics and its Applications 71

[3] Johann Cigler, Hankel determinants, Catalan numbers and Fibonacci polynomials, arXiv:1801.05608

[4] Johann Cigler, Shifted Hankel determinants of Catalan numbers and related results II: Backward shifts, arxiv: 2306.07733

[5] Johann Cigler, Some experimental observations about Hankel determinants of convolution powers of Catalan numbers, arXiv: 2308.07642

[6] Johann Cigler and Christian Krattenthaler, Some determinants of path generating functions, Adv. Appl. Math. 46 (2011), 144-174

[7] Christian Krattenthaler, Advanced determinant calculus: A complement, Linear Algebra Appl. 411 (2005), 68-166

[8] OEIS, The On-Line Encyclopedia of Integer Sequences,